\documentclass[11pt]{article}
\usepackage{amssymb}
\usepackage{amsmath}
\usepackage{amsfonts}
\usepackage{bbm}
\usepackage{enumerate}
\newtheorem{theorem}{Theorem}[section]

\newtheorem{remark}[theorem]{Remark}

\begin{document}

\title{Criteria of valid line sum arrays for multidimensional matrices}

\author{Hyun Kwang Kim\thanks{The work of Hyun Kwang Kim was supported by Basic Research Program through the National Research Foundation of Korea (NRF) funded by the Ministry of Education, Science and Technology (grant \# 2009-0089826).}\\
Department of Mathematics\\ 
POSTECH\\
Pohang 790-784, Korea\\
hkkim@postech.ac.kr
\and
Joon Yop Lee\thanks{Corresponding author}\\
Department of Mathematics\\
POSTECH\\
Pohang 790-784, Korea\\
flutelee@postech.ac.kr}

\maketitle

\begin{abstract}
Gale and Ryser found a criterion for existence of a binary matrix with given row and column sums.
Mirsky extended the theorem of Gale and Ryser to $q$-ary matrices.
In this paper, we are interested in higher dimensional extension of these theorems.
We first introduce multidimensional matrices as a higher dimensional generalization of matrices.
We next replace the concept of row and columns in matrices by lines in multidimensional matrices.
We finally find a criterion for existence of a $q$-ary multidimensional matrices with given line sums.

\end{abstract}

\section{Introduction}
A matrix is called \emph{binary} if each of its entries is either $0$ or $1$.
Let $M=(m_{ij})$ be an $n_1 \times n_2$ binary matrix. 
The row sum vector of $M$ is a vector $\bf{r}$ = $(r_1,r_2, \ldots,r_{n_1})$ of length $n_1$ such that $r_i=\sum_{j=1}^{n_2}m_{ij}$. 
Similarly, the column sum vector of $M$ is a vector $\bf{c}$ = $(c_1,c_2, \ldots, c_{n_2})$ of length $n_2$ such that $c_j=\sum_{i=1}^{n_1}m_{ij}$.
A pair $(\bf{u},\bf{v})$, where $\bf{u}$ is a vector of length $n_1$ and $\bf{v}$ is a vector of length $n_2$ is called \emph{valid} if there is an $n_1 \times n_2$ binary matrix $M$ whose row (resp. column) sum vector is $\bf{u}$ (resp. $\bf{v}$).
One of the main problem on row and column sums of matrices is to find a criterion for a pair $(\bf{u},\bf{v})$ to be valid.
The following theorem due to Gale \cite{gale} and Ryser \cite{ryser} 
is a fundamental result in this line of research.

\begin{theorem}\label{theorem1.1}
Let $r_1$, $r_2$,\,\ldots, $r_{n_1}$ and $c_1$, $c_2$,\,\ldots, $c_{n_2}$ be nonnegative integers such that
$$\begin{cases}
0\leq r_{i_1}\leq n_2, 0\leq c_{i_2}\leq n_1\\
\underset{i_1=1}{\overset{n_1}{\sum}}r_{i_1}=\underset{i_2=1}{\overset{n_2}{\sum}}c_{i_2}\\
c_1\ge c_2\ge\cdots\ge c_{n_2}
\end{cases}.$$
Let $\overline{M}$ be a binary matrix of size $n_1\times n_2$ with the row sums $r_1$, $r_2$, \ldots, $r_{n_1}$ such that
$$\overline{M}(i_1,i_2)\ge\overline{M}(i_1,i_2+1)$$
If the column sums of $\overline{M}$ are $\overline{c}_1$, $\overline{c}_2$, \ldots, $\overline{c}_{n_2}$, then there exists a binary matrix $M$ of size $n_1\times n_2$ with the row sums $r_1$, $r_2$,\,\ldots, $r_{n_1}$ and column sums $c_1$, $c_2$,\,\ldots, $c_{n_2}$ if and only if $$c_1+c_2+\cdots+c_j\leq \overline{c}_1+\overline{c}_2+\cdots+\overline{c}_j\ (1\leq j\leq n_1).$$
\end{theorem}

A \emph{$q$-ary matrix} is a matrix each of whose entries is in $\{0,1,\ldots,q-1\}$.
Mirsky \cite{mirsky,mirsky1} extended the above theorem of Gale and Ryser to $q$-ary matrices, that is, he found a criterion for existence of a $q$-ary matrix with given row and column sums.

The goal of this paper is to consider a higher dimensional extension of the problems described above.
The position of each entry in a matrix is determined by its row and column positions.
Thus we may regard matrices as $2$-dimensional objects. Using this observation, we can generalize matrices as follows:
For positive integers $n_1,n_2,\ldots,n_d$ let $M=[M(i_1,i_2,\ldots,i_d)]$ be a map defined by
$$M:\prod_{i=1}^d\{1,2,\ldots,n_i\}\rightarrow\{0,1,2,\ldots,q-1\}.$$
We call such an $M$ a \emph{multidimensional matrix}. Multidimensional matrices are higher dimensional generalizations of matrices.
So we naturally wonder occurring phenomena when we extend row sum and column sum problems on matrices to multidimensional matrices.
By employing concepts in multidimensional matrices, we can generalize rows and columns in matrices to \emph{lines} in multidimensional matrices.
(The exact definition of line will be given in the next section.)
%These generalizations suggest us to consider several criteria on multidimensional matrices, which are similar to the case of ordinary matrices.
Haber \cite{haber, haber1} and Ryser \cite{ryser, ryser1, ryser2} gave a criterion for existence of a binary matrix with given row and column sums by constructing a binary matrix with given row and column sums. We generalize the methods of Haber and Ryser and establish several criteria for existence of a $q$-ary multidimensional matrix with given line sums by a constructive method.

%Main method to achieve our goal is an extension of Haber's and Ryser's one, that is, we construct a multidimensional matrix with given line sums.

\section{Preliminaries\label{ldefinition}}
This section is devoted to generalization of various concepts in matrices to those in multidimensional matrices.
This generalization gives multidimensional matrix analogs of corresponding concepts in matrices.
We also provide the main idea to construct a multidimensional matrices with given line sums.

\subsection{Multidimensional matrices}
Let $n$ be a positive integer. We denote by $[n]$ the set $\{1,2,\ldots,n\}$. For a $d$-dimensional vector $\mathbf{n}_d=(n_1,n_2,\ldots,n_d)$ with positive integer entries, we set $[\mathbf{n}_d]=\prod_{i=1}^d[n_i]$ and denote an element of $[\mathbf{n}_d]$ by $\mathbf{i}_d=(i_1,i_2, \ldots ,i_d)$.
Fixing $d-1$ entries $i_1$,\,\ldots,$i_{j-1}$, $i_{j+1}$,\,\ldots, $i_d$ among $d$ coordinates in $\mathbf{i}_d$, we define a \emph{line} of $[\mathbf{n}_d]$ to be a subset
$$\mathbf{L}_j(\mathbf{i}_d)=\{\mathbf{i}_d\mid i_j\in[n_j]\}.$$
We call it the \emph{$j$th line} of $[\mathbf{n}_d]$ at $\mathbf{i}_d$. Note that each $j$th line $(1\leq j\leq d)$ is invariant for any value of $i_j$. However, to simplify the notation, we denote a line by $\mathbf{L}_j(\mathbf{i}_d)$ by abusing notation.
We define $\mathcal{L}(\mathbf{n}_d)$ to be the set of lines in $[\mathbf{n}_d]$.

As we have mentioned in introduction, we define a \emph{$q$-ary} (\emph{$d$-dimensional}) \emph{matrix} of size $n_1\times n_2\times \cdots\times n_d$ to be a map $M:[\mathbf{n}_d]\rightarrow \{0,1,\dots,q-1\}$ and denote it by
$$M=[M(\mathbf{i}_d)\mid\mathbf{i}_d\in[\mathbf{n}_d]].$$
We simply call $M$ an \emph{$(\mathbf{n}_d,q)$-matrix} or a \emph{$d$-matrix}, and we assume that matrices mean multidimensional matrices.
For a subset $\mathbf{N}$ of $[\mathbf{n}_d]$, we call $M(\mathbf{N})=[M(\mathbf{i}_d)\mid\mathbf{i}_d\in\mathbf{N}]$ a \emph{subarray} of $M$.
Along the same line as the case of $[\mathbf{n}_d]$, we define the \emph{$j$th line} of $M$ at $\mathbf{i}_d$ to be a subarray $M(\mathbf{L}_j(\mathbf{i}_d))$.
We denote $|M(\mathbf{L}_j(\mathbf{i}_d))|=\sum_{\mathbf{x}_d\in\mathbf{L}_j}|M(\mathbf{x}_d)|$ and call it a \emph{line sum}.

To define a concept relevant to the line sum set of a matrix, we let the \emph{line sum array} of $M$ be
$$S(M)=[|M(\mathbf{L})|\mid\mathbf{L}\in\mathcal{L}(\mathbf{n}_d)].$$
Note that for each $1\leq j\leq d$ the array $[|M(\mathbf{L}_j(\mathbf{i}_d))|\mid i_k\in[n_k]\ (k\neq j)]$ is a matrix of size $n_1\times \cdots n_{j-1}\times n_{j+1}\times\cdots\times n_d$.
Thus a line sum array is composed of matrices.

\begin{remark}\label{remark}

The definitions of line and line sum imply the followings:
Let $M$ be a $q$-ary matrix of size $n_1\times n_2$.
\begin{enumerate}[1.]
\item The first line of $M$ at $\mathbf{i}_2$ is the $i_2$th column of $M$ and the second line of $M$ at $\mathbf{i}_2$ is the $i_1$th row of $M$.

\item The line sums of $M$ satisfy
$$\sum_{i_2=1}^{n_2}|M(\mathbf{L}_1(\mathbf{i}_2))|=\sum_{i_1=1}^{n_1}|M(\mathbf{L}_2(\mathbf{i}_2))|.$$
Thus, if $M^{\prime}$ is a $q$-ary matrix of size $n_1\times n_2\times\cdots\times n_d$, then for each $\mathbf{i}_d\in[\mathbf{n}_d]$ and $j,j^{\prime}\in[d]$ the line sums of $M^{\prime}$ satisfies
$$\sum_{i_{j^{\prime}}=1}^{n_{j^{\prime}}}|M^{\prime}(\mathbf{L}_j(\mathbf{i}_d))|=\sum_{i_j=1}^{n_j}|M^{\prime}(\mathbf{L}_{j^{\prime}}(\mathbf{i}_d))|.$$
\end{enumerate}
\end{remark}

\subsection{Line sum arrays}

A \emph{$q$-ary} (\emph{$d$-dimensional}) \emph{line sum array} of size $n_1\times n_2\times\cdots\times n_d$ is an array $S=[S(\mathbf{L})\mid\mathbf{L}\in\mathcal{L}(\mathbf{n}_d)]$ with nonnegative integer entries such that
\begin{enumerate}[1.]
\item\label{linesumcondition1}
$0\leq S(\mathbf{L})\leq |\mathbf{L}|\cdot(q-1)$,

\item\label{linesumcondition2}
$\sum_{i_{j^{\prime}}=1}^{n_{j^{\prime}}}S(\mathbf{L}_j(\mathbf{i}_d))
=\sum_{i_j=1}^{n_j}S(\mathbf{L}_{j^{\prime}}(\mathbf{i}_d))$ $(\mathbf{i}_d\in[\mathbf{n}_d],j,j^{\prime}\in[d])$.
\end{enumerate}
The concept of line sum array corresponds to that of the line sum array of a matrix.
Thus we also call $S$ an \emph{$(\mathbf{n}_d,q)$-line sum array} or a \emph{$d$-line sum array}.
We say that $S$ is \emph{valid} if there is a $q$-ary matrix $M$ of size $n_1\times n_2\times\cdots\times n_d$ satisfying $$S(M)=S.$$
We call $M$ a \emph{matrix of $S$} and denote $M=M(S)$.
Note that the conditions in the definition of line sum array are simple necessary conditions for a line sum array to be valid.
(For condition \ref{linesumcondition1} refer to the definition of $q$-ary matrix and for condition \ref{linesumcondition2} refer to Remark \ref{remark}.)

Let $S$ be a $q$-ary line sum array of size $n_1\times n_2\times\cdots\times n_d$.
There is a unique $q$-ary matrix $\overline{M}$ of size $n_1\times n_2\times\cdots\times n_d$ that satisfies the following conditions:
\begin{enumerate}[1.]
\item
$|\overline{M}(\mathbf{L}_d(\mathbf{i}_d))|=S(\mathbf{L}_d(\mathbf{i}_d))$ for each $d$th line $\overline{M}(\mathbf{L}_d(\mathbf{i}_d))$ of $\overline{M}$.

\item
$\overline{M}(\mathbf{i}_{d-1},i_d)\ge\overline{M}(\mathbf{i}_{d-1},i_d+1)$.

\item
Each $d$th line $\overline{M}(\mathbf{L}_d(\mathbf{i}_d))$ contains
at most one entry whose value is neither $0$ nor $q-1$.
\end{enumerate}
We call such a matrix a \emph{maximal matrix} and denote $\overline{M}=\overline{M}(S)$ and $\overline{S}=S(\overline{M})$.
Figure \ref{tmaximal} is an example of a maximal matrix.
\begin{figure}[h]
$$\begin{pmatrix}
0&0&0&0&0&0&0&0&0&0&0\\
1&0&0&0&0&0&0&0&0&0&0\\
2&0&0&0&0&0&0&0&0&0&0\\
2&2&2&2&2&2&0&0&0&0&0\\
2&2&2&2&2&2&2&2&2&2&0\\
2&2&2&2&1&0&0&0&0&0&0\\
2&2&2&2&2&2&2&2&2&2&1\\
2&2&2&2&2&2&2&2&2&2&2\\
2&2&2&1&0&0&0&0&0&0&0\\
2&2&2&2&2&2&2&2&1&0&0
\end{pmatrix}$$
\caption{A $3$-ary maximal matrix of size $10\times 11$\label{tmaximal}}
\end{figure}

\subsection{The main idea of theorems}
Let $S$ be a $q$-ary line sum array of size $n_1\times n_2$. For a permutation $\sigma$ of $[n_2]$, we let $S^{\sigma}$ be a $q$-ary line sum array of size $n_1\times n_2$ defined by
$$\begin{cases}
S^{\sigma}(\mathbf{L}_1(\mathbf{i}_2))=S(\mathbf{L}_1(i_1,\sigma(i_2)))\\
S^{\sigma}(\mathbf{L}_2(\mathbf{i}_2))=S(\mathbf{L}_2(\mathbf{i}_2))
\end{cases}\ (\mathbf{i}_2\in[\mathbf{n}_2]).$$
Suppose that $S$ is valid.
The definition of maximal matrix implies that
\begin{align*}
\sum_{i_2=1}^{a_2}S(\mathbf{L}_1(\mathbf{i}_2))&=
\sum_{i_2=1}^{a_2}\sum_{i_1=1}^{n_1}M(S)(\mathbf{i}_2)\\
&\leq\sum_{i_2=1}^{a_2}\sum_{i_1=1}^{n_1}\overline{M}(S)(\mathbf{i}_2)
=\sum_{i_2=1}^{a_2}\overline{S}(\mathbf{L}_1(\mathbf{i}_2))\ (1\leq a_2\leq n_2).
\end{align*}
Therefore
\begin{equation*}
\sum_{i_2=1}^{a_2}S(\mathbf{L}_1(\mathbf{i}_2))\leq
\sum_{i_2=1}^{a_2}\overline{S}(\mathbf{L}_1(\mathbf{i}_2))\ (1\leq a_2\leq n_2).
\end{equation*}
Similarly, for each permutation $\sigma$ of $[n_2]$ the line sum array $S$ satisfies
\begin{equation}\label{lcondition}
\sum_{i_2=1}^{a_2}S^{\sigma}(\mathbf{L}_1(\mathbf{i}_2))\leq
\sum_{i_2=1}^{a_2}\overline{S^{\sigma}}(\mathbf{L}_1(\mathbf{i}_2))\ (1\leq a_2\leq n_2).
\end{equation}
We say that a $q$-ary line sum array $S$ of size $n_1\times n_2$ is \emph{compatible} if $S$ satisfies inequalities (\ref{lcondition}) for each permutation $\sigma$ of $[n_2]$.
Using the definition of compatibility, we can restate Theorem \ref{theorem1.1} in the following way.
\begin{theorem}
A binary $2$-line sum array $S$ is valid if and only if $S$ is compatible.
\end{theorem}

Let $S$ be a $q$-ary line sum array of size $n_1\times n_2\times\cdots\times n_d$.
We define a \emph{$2$-line sum subarray} of $S$ to be an array $S^{\prime}$ such that $$S^{\prime}=[S(\mathbf{L}_j(\mathbf{i}_d))\mid j\in\{j_1,j_2\},i_{j^{\prime}  }\in[n_{j^{\prime}}]\ \mbox{for}\ j^{\prime}\in[d]\setminus\{j_1,j_2\}]$$
for some distinct $j_1$ and $j_2$.
Note that $S^{\prime}$ corresponds to the line sum array of a $2$-matrix of size $n_{j_1}\times n_{j_2}$.
We say that $S$ is \emph{compatible} if each of its $2$-line sum subarray is compatible.
We show that compatibility is a necessary and sufficient condition for $S$ to be valid, which is a generalization of Theorem \ref{theorem1.1}.
Since a valid line sum array is trivially compatible, in the proof of each theorem we omit to prove the trivial fact that ``a valid line sum array is compatible.''

\section{Criteria for $2$-matrices\label{lsection}}

The main method that we will use  is based on Haber's \cite{haber, haber1} and Ryser's \cite{ryser, ryser1, ryser2} constructive proof of Theorem \ref{theorem1.1}.
Therefore we review Haber's and Ryser's proof by interpreting their proof in terms of our new terminologies.
Using this new interpretation of Haber's and Ryser's proof, we establish several criteria for existence of $q$-ary multidimensional matrices.

\subsection{Binary $2$-matrices}
Haber \cite{haber, haber1} and Ryser \cite{ryser, ryser1, ryser2} proved Theorem \ref{theorem1.1} by constructing a binary $2$-matrix with a given $2$-line sum array.
Our method to derive criteria for multidimensional matrices is an extension of Haber's and Ryser's proof, that is, we construct a $q$-ary multidimensional matrix with a given line sum array.
Thus we briefly introduce here a proof translated from Haber's and Ryser's proof by borrowing terminologies defined in Section \ref{ldefinition}.

\begin{pftheorem}
Let $S$ be a compatible binary line sum array of size $n_1\times n_2$.
To show that $S$ is valid, we construct a binary matrix $M$ of size $n_1\times n_2$ with $S(M)=S$ by using induction on the number $n_2$.
If $n_2=1$, then $S$ is trivially valid.
Thus we assume that $n_2\ge 2$.

The idea of our construction is the following:
Let $\overline{M}$ be the maximal matrix of $S$.
A (\emph{$2$-dimensional}) \emph{shift operation} is successive shifts of the rightmost one in a row of $\overline{M}$ to the last column in the same row.
We construct $M$ from $\overline{M}$ by applying shift operations to $\overline{M}$.
Since we will frequently use operations similar to shift operation to prove theorems for multidimensional matrices, we formally define a ($2$-dimensional) shift operation to be successive shifts of an entry $\overline{M}(\mathbf{i}_2)$ in $\overline{M}$ with
$$\overline{M}(\mathbf{i}_2)=1\ \text{and}\ \overline{M}(i_1,i_2+1)=0$$
along the 2nd line $\overline{M}(\mathbf{L}_2(\mathbf{i}_2))$ that yields $M(i_1,n_2)=1$. Figure \ref{shift1} is an example of shift operation.
\begin{figure}[h]
$$\begin{pmatrix}
1&1&1&0&0&0&0&0&0\\
1&1&1&1&1&\mathbf{1}&0&0&0\\
1&1&1&1&1&1&0&0&0\\
1&0&0&0&0&0&0&0&0
\end{pmatrix}
\longrightarrow
\begin{pmatrix}
1&1&1&0&0&0&0&0&0\\
1&1&1&1&1&0&0&0&\mathbf{1}\\
1&1&1&1&1&1&0&0&0\\
1&0&0&0&0&0&0&0&0
\end{pmatrix}
$$
\caption{A shift operation\label{shift1}}
\end{figure}

We apply $S(\mathbf{L}_1(i_1,n_2))$ shift operations to $\overline{M}(S)$ according to the following order:
We apply shift operations from the rightmost one.
And if two ones are located in the same column, then we apply shift operations from the upper one.

Let $\mathbf{I}$ be the set of indexes $i_1$ of rows at which shift operations occur.
It follows from our construction that $\mathbf{I}$ satisfies the following conditions:
\begin{enumerate}[{1.}]
\item
$|\mathbf{I}|=S(\mathbf{L}_1(i_1,n_2))$.

\item
Any two $i_1\in\mathbf{I}$ and $i^{\prime}_1\in [n_1]\setminus \mathbf{I}$ satisfies one of the followings;
\begin{enumerate}[1.]
\item
$S(\mathbf{L}_2(\mathbf{i}_2))>S(\mathbf{L}_2(i^{\prime}_1,i_2))$.

\item
$S(\mathbf{L}_2(\mathbf{i}_2))=S(\mathbf{L}_2(i^{\prime}_1,i_2))$ and $i_1<i^{\prime}_1$.
\end{enumerate}
\end{enumerate}
Shift operations on the rows of $\overline{M}(S)$ whose indexes are in the set $\mathbf{I}$ yields a binary matrix $M^{\prime}=[\overline{M}_{n_2-1},M_2]$ of size $n_1\times n_2$ with line sums
$$\begin{cases}
|M^{\prime}(\mathbf{L}_1(i_1,n_2))|=S(\mathbf{L}_1(i_1,n_2))\\
|M^{\prime}(\mathbf{L}_2(\mathbf{i}_2))|=S(\mathbf{L}_2(\mathbf{i}_2))
\end{cases}\ (1\leq i_1\leq n_1)$$
such that $\overline{M}_{n_2-1}$ is a binary maximal matrix of size $n_1\times (n_2-1)$ and $M_1$ is a binary matrix of size $n_1\times 1$.
To show that this construction of $M^{\prime}$ is valid, we need to verify that the number $S(\mathbf{L}_1(i_1,n_2))$ satisfies
\begin{equation}\label{compatible1}
\overline{S}(\mathbf{L}_1(i_1,n_2))\leq S(\mathbf{L}_1(i_1,n_2))\leq \overline{S}(\mathbf{L}_1(i_1,1)).
\end{equation}

The compatibility of $S$ implies that
$$\begin{cases}
\sum_{i_2=1}^{n_2-1}S(\mathbf{L}_1(\mathbf{i}_2))\leq
\sum_{i_2=1}^{n_2-1}\overline{S}(\mathbf{L}_1(\mathbf{i}_2))\\
\sum_{i_2=1}^{n_2}S(\mathbf{L}_1(\mathbf{i}_2))=
\sum_{i_2=1}^{n_2}\overline{S}(\mathbf{L}_1(\mathbf{i}_2))
\end{cases}.$$
Thus
$$\overline{S}(\mathbf{L}_1(i_1,n_2))\leq S(\mathbf{L}_1(i_1,n_2)).$$
Let $\sigma$ be a permutation of $[n_2]$ obeying $\sigma(1)=n_2$.
Then the line sum $S(\mathbf{L}_1(i_1,n_2))$ satisfies
$$S(\mathbf{L}_1(i_1,n_2))=
S^{\sigma}(\mathbf{L}_1(i_1,1))\leq
\overline{S^{\sigma}}(\mathbf{L}_1(i_1,1))
=\overline{S}(\mathbf{L}_1(i_1,1)),$$
which shows that
$$S(\mathbf{L}_1(i_1,n_2))\leq \overline{S}(\mathbf{L}_1(i_1,1)).$$
Therefore $S(\mathbf{L}_1(i_1,n_2))$ satisfies (\ref{compatible1}).

The line sum array $S_{n_2-1}$ of size $n_1\times (n_2-1)$ defined by
$$\begin{cases}
S_{n_2-1}(\mathbf{L}_1(\mathbf{i}_2))=S(\mathbf{L}_1(\mathbf{i}_2))&(1\leq i_2\leq n_2-1)\\
S_{n_2-1}(\mathbf{L}_2(\mathbf{i}_2))=|\overline{M}_{n_2-1}(\mathbf{L}_2(\mathbf{i}_2))|&(1\leq i_1\leq n_1)\end{cases}$$
is compatible.
Therefore induction on the number $n_2$ gives a binary matrix $M_{n_2-1}$ of size $n_1\times (n_2-1)$ such that the binary matrix $M=[M_{n_2-1},M_1]$ of size $n_1\times n_2$ satisfies
$$S(M)=S.$$
This proves the theorem.
\end{pftheorem}

\subsection{$q$-ary $2$-matrices}

A $q$-ary generalization of Theorem \ref{theorem1.1} is the following Mirsky's theorem \cite{mirsky}.
\begin{theorem}\label{qlmatrix}
Let $S$ be a $q$-ary $2$-line sum array.
Then $S$ is valid if and only if $S$ is compatible.
\end{theorem}
Dias da Silva et al. proved Theorem \ref{qlmatrix} by constructing a $q$-ary $2$-matrix with a given line sum array \cite{dias}.
Since we consider a generalization of this theorem to multidimensional matrices, we provide here a proof of Theorem \ref{qlmatrix} that employs the concept of multidimensional matrix.

For this, we introduce the following representation of $q$-ary matrices in terms of binary matrices:
Let $M$ be a $q$-ary matrix of size $n_1\times n_2\times\cdots\times n_d$.
We define the \emph{binary representation of $M$} to be a binary matrix $M^b$ of size $n_1\times n_2\times\cdots\times n_d\times(q-1)$ such that
\begin{equation}\label{qary}\begin{cases}
M^b(\mathbf{i}_d,i_{d+1})=1&(1\leq i_{d+1}\leq M(\mathbf{i}_d))\\
M^b(\mathbf{i}_d,i_{d+1})=0&(M(\mathbf{i}_d)+1\leq i_{d+1}\leq q-1)
\end{cases}.\end{equation}

\begin{pfqlmatrix}
Let $S$ be a compatible $q$-ary line sum array of size $n_1\times n_2$.
To show that $S$ is valid, we construct a $q$-ary matrix $M$ of size $n_1\times n_2$ with $S(M)=S$ by using induction on the number $n_2$.
If $n_2=1$, then $S$ is trivially valid.
Thus we assume that $n_2\ge 2$.

Similar to the binary case, the idea of construction is the following: Let $\overline{M}$ be the maximal matrix of $M$.
A (\emph{$q$-ary $2$-dimensional}) \emph{shift operation} is successive shifts of an entry $\overline{M}^b(\mathbf{i}_3)$ of $\overline{M}^b$ satisfying
$$\overline{M}^b(\mathbf{i}_3)=1\ \text{and}\ \overline{M}^b(i_1,i_2+1,i_3)=0$$
along the $2$nd line $\overline{M}^b(\mathbf{L}_2(\mathbf{i}_3))$ that yields $M^b(i_1,n_2,i_3)=1$.
Figure \ref{shift2} is an example of a $3$-ary shift operation.
\begin{figure}[h]
$$\begin{pmatrix}
2&2&1&0&0&0&0&0&0\\
2&2&2&2&2&\mathbf{2}&0&0&0\\
2&2&2&2&2&1&0&0&0\\
2&0&0&0&0&0&0&0&0
\end{pmatrix}
\longrightarrow
\begin{pmatrix}
2&2&1&0&0&0&0&0&0\\
2&2&2&2&2&\mathbf{1}&0&0&\mathbf{1}\\
2&2&2&2&2&1&0&0&0\\
2&0&0&0&0&0&0&0&0
\end{pmatrix}
$$
\caption{A $3$-ary shift operation\label{shift2}}
\end{figure}
We construct $M$ from $\overline{M}$ by applying shift operations to $\overline{M}$.

We apply shift operations to $\overline{M}(S)^b$ to construct a $q$-ary matrix $M^{\prime}=[\overline{M}_{n_2-1},M_1]$ of size $n_1\times n_2$ with line sums
$$\begin{cases}
|M^{\prime}(\mathbf{L}_1(i_1,n_2))|=S(\mathbf{L}_1(i_1,n_2))\\
|M^{\prime}(\mathbf{L}_2(\mathbf{i}_2))|=S(\mathbf{L}_2(\mathbf{i}_2))
\end{cases}\ (1\leq i_1\leq n_1)$$
such that $\overline{M}_{n_2-1}$ is a $q$-ary maximal matrix of size $n_1\times (n_2-1)$ and $M_1$ is a $q$-ary matrix of size $n_1\times 1$.
Since we can only apply shift operations to $S(\mathbf{L}_1(i_1,n_2))$ $2$nd lines of $\overline{M}(S)^b$, we need to choose $S(\mathbf{L}_1(i_1,n_2))$ $2$nd lines of $\overline{M}(S)^b$.

For this, we define $\mathbf{I}$ to be the set of $(i_1,i_3)\in[n_1]\times[q-1]$ that satisfies the following conditions:
\begin{enumerate}[{1.}]
\item
$|\mathbf{I}|=S(\mathbf{L}_1(i_1,n_2,i_3))$.

\item
Any two $(i_1,i_3)\in\mathbf{I}$ and $(i^{\prime}_1,i^{\prime}_3)\in[n_1]\times [q-1]\setminus \mathbf{I}$ satisfies one of the followings:
\begin{enumerate}[(1)]
\item
$\overline{M}(S)^b(\mathbf{L}_2(\mathbf{i}_3))|>
|\overline{M}(S)^b(\mathbf{L}_2(i^{\prime}_1,i_2,i^{\prime}_3))|$.

\item
$|\overline{M}(S)^b(\mathbf{L}_2(\mathbf{i}_3))|=
|\overline{M}(S)^b(\mathbf{L}_2(i^{\prime}_1,i_2,i^{\prime}_3))|$ and $i_3<i^{\prime}_3$.

\item
$|\overline{M}(S)^b(\mathbf{L}_2(\mathbf{i}_3))|=
|\overline{M}(S)^b(\mathbf{L}_2(i^{\prime}_1,i_2,i^{\prime}_3))|$, $i_3=i^{\prime}_3$, and $i_1<i^{\prime}_1$.
\end{enumerate}
\end{enumerate}
The definition of $\mathbf{I}$ implies that applying shift operations to the lines $\overline{M}(S)^b(\mathbf{L}_2(\mathbf{i}_3))$ for $(i_1,i_3)\in\mathbf{I}$ yield $M^{\prime}$.
To show that this construction is valid, we need to verify that the number $S(\mathbf{L}_1(i_1,n_2))$ satisfies
\begin{equation}\label{compatible2}
\overline{S}(\mathbf{L}_1(i_1,n_2))\leq S(\mathbf{L}_1(i_1,n_2))\leq\overline{S}(\mathbf{L}_1(i_1,1)).
\end{equation}

The compatibility of $S$ implies that
$$\begin{cases}
\sum_{i_2=1}^{n_2-1}S(\mathbf{L}_1(\mathbf{i}_2))\leq
\sum_{i_2=1}^{n_2-1}\overline{S}(\mathbf{L}_1(\mathbf{i}_2))\\
\sum_{i_2=1}^{n_2}S(\mathbf{L}_1(\mathbf{i}_2))=
\sum_{i_2=1}^{n_2}\overline{S}(\mathbf{L}_1(\mathbf{i}_2))
\end{cases}.$$
Thus we obtain
$$\overline{S}(\mathbf{L}_1(i_1,n_2))\leq S(\mathbf{L}_1(i_1,n_2)).$$
Let $\sigma$ be a permutation of $[n_2]$ obeying $\sigma(1)=n_2$.
Then the line sum $S(\mathbf{L}_1(i_1,n_2))$ satisfies
$$S(\mathbf{L}_1(i_1,n_2))=S^{\sigma}(\mathbf{L}_1(i_1,1))\leq
\overline{S^{\sigma}}(\mathbf{L}_1(i_1,1))=\overline{S}(\mathbf{L}_1(i_1,1)),$$
which shows that
$$S(\mathbf{L}_1(i_1,n_2))\leq \overline{S}(\mathbf{L}_1(i_1,1)).$$
Therefore $S(\mathbf{L}_1(i_1,n_2))$ satisfies (\ref{compatible2}).

The $q$-ary line sum array $S_{n_2-1}$ of size $n_1\times (n_2-1)$ defined by
$$\begin{cases}
S_{n_2-1}(\mathbf{L}_1(\mathbf{i}_2))=S(\mathbf{L}_1(\mathbf{i}_2))&(1\leq i_2\leq n_2-1)\\
S_{n_2-1}(\mathbf{L}_2(\mathbf{i}_2))=|\overline{M}_{n_2-1}(\mathbf{L}_2(\mathbf{i}_2))|&(1\leq i_1\leq n_1)\\
\end{cases}$$
is compatible.
Therefore induction on the number $n_2$ gives a $q$-ary matrix $M_{n_2-1}$ of size $n_1\times (n_2-1)$ such that the $q$-ary matrix $M=[M_{n_2-1},M_1]$ of size $n_1\times n_2$ satisfies
$$S(M)=S.$$
This proves the theorem.
\end{pfqlmatrix}

\section{Criteria for multidimensional matrices}
We have defined multidimensional matrix by generalizing the concept of $2$-matrix and constructed a $q$-ary $2$-matrix with a given line sum array by using the concept of binary representation.
On the basis of these definition and construction, we establish multidimensional generalizations of Theorems \ref{theorem1.1} and \ref{qlmatrix}.
We also provide such theorems for symmetric multidimensional matrices.

\medskip
The following theorem gives a criterion for existence of a matrix with a given line sum array.
\begin{theorem}\label{qmlmatrix}
Let $S$ be a line sum array.
Then $S$ is valid if and only if $S$ is compatible.
\end{theorem}
The main idea to prove Theorem \ref{qmlmatrix} is similar to the case of $2$-matrices, that is, we use a multidimensional version of shift operation.
It is defined as follows:
A (\emph{$q$-ary $d$-dimensional}) \emph{shift operation} is successive shifts of an entry $\overline{M}(S)^b(\mathbf{i}_{d+1})$ of $\overline{M}(S)^b$ satisfying
$$\overline{M}(S)^b(\mathbf{i}_{d+1})=1\ \text{and}\
\overline{M}(S)^b(\mathbf{i}_{d-1},i_d+1,i_{d+1})=0$$
along a $d$th line $\overline{M}(S)^b(\mathbf{L}_d(\mathbf{i}_{d+1}))$ that yields $M^b(\mathbf{i}_{d-1},n_d,i_{d+1})=1$.
We construct an $(\mathbf{n}_d,q)$-matrix $M$ with $S(M)=S$ from $\overline{M}(S)^b$ by applying shift operations to $\overline{M}(S)^b$.

\begin{pfqmlmatrix}
Let $S$ be a compatible $(\mathbf{n}_d,q)$-line sum array.
We construct an $(\mathbf{n}_d,q)$-matrix $M$ with $S(M)=S$ by using double induction on the numbers $d$ and $n_d$.

Suppose that $d=2$.
By Theorem \ref{qlmatrix}, $S$ is valid.
From now on we assume that $d\ge 3$.

Suppose that $n_d=1$.
By induction on the number $d$, we can construct an $(\mathbf{n}_{d-1},q)$-matrix $M_{d-1}$ such that
$$|M_{d-1}(\mathbf{L}_j(\mathbf{i}_{d-1}))|=S(\mathbf{L}_j(\mathbf{i}_{d-1},1))\ (1\leq j\leq d-1,\mathbf{i}_{d-1}\in[\mathbf{n}_{d-1}]).$$
Therefore the matrix $M$ defined by $M(\mathbf{i}_{d-1},1)=M_{d-1}(\mathbf{i}_{d-1})$ $(\mathbf{i}_{d-1}\in[\mathbf{n}_{d-1}])$ satisfies $S(M)=S$, which shows that $S$ is valid.

Suppose that $n_d\ge 2$.
We apply shift operations to the matrix $\overline{M}(S)^b$ to construct an $(\mathbf{n},q)$-matrix $M^{\prime}=[\overline{M}_{n_d-1},M_1]$ with line sums
$$\begin{cases}
|M^{\prime}(\mathbf{L}_d(\mathbf{i}_d))|
=S(\mathbf{L}_d(\mathbf{i}_d))&(\mathbf{i}_d\in[\mathbf{n}_d])\\
|M^{\prime}(\mathbf{L}_j(\mathbf{i}_{d-1},n_d))|
=S(\mathbf{L}_j(\mathbf{i}_{d-1},n_d))&(1\leq j\leq d-1,\mathbf{i}_{d-1}\in[\mathbf{n}_{d-1}])
\end{cases}$$
such that $\overline{M}_{n_d-1}$ is a maximal $((\mathbf{n}_{d-1},n_d-1),q)$-matrix and $M_1$ is an $((\mathbf{n}_{d-1}, 1),q)$-matrix.
Similar to the case of $2$-matrices, for each $1\leq j\leq d-1$ and $\mathbf{i}_{d-1}\in[\mathbf{n}_{d-1}]$ we can apply shift operations to only $S(\mathbf{L}_j(\mathbf{i}_{d-1},n_d))$ $d$th lines of $\overline{M}(S)^b$ among the lines
$\overline{M}(\mathbf{L}_d(\mathbf{i}_{d+1}))$ $(1\leq i_j\leq n_j)$.
Thus for each $1\leq j\leq d-1$ and $\mathbf{i}_{d-1}\in[\mathbf{n}_{d-1}]$ we need to choose $S(\mathbf{L}_j(\mathbf{i}_{d-1},n_d))$ $d$th lines of $\overline{M}(S)^b$ among the lines
$\overline{M}(\mathbf{L}_d(\mathbf{i}_{d+1}))$ $(1\leq i_j\leq n_j)$.

For this construction, we define $\mathbf{I}$ to be the set of $(\mathbf{i}_{d-1},i_{d+1})\in[\mathbf{n}_{d-1}]\times[q-1]$ that satisfies the following conditions:
\begin{enumerate}[1.]
\item
For each $1\leq j\leq d-1 $ and $\mathbf{i}_{d-1}\in[\mathbf{n}_{d-1}]$, the number of $(\mathbf{i}_{d-1},i_{d+1})$ in $\mathbf{I}$ is $S(\mathbf{L}_j(\mathbf{i}_{d-1},n_d))$.

\item
Any two $(\mathbf{i}_{d-1},i_{d+1})\in\mathbf{I}$ and $(\mathbf{i}^{\prime}_{d-1},i^{\prime}_{d+1})\in[\mathbf{n}_{d-1}]\times [q-1]\setminus \mathbf{I}$ satisfy one of the followings:
\begin{enumerate}[(1)]
\item
$|\overline{M}(S)^b(\mathbf{L}_d(\mathbf{i}_{d+1}))|>|\overline{M}(S)^b(\mathbf{L}_d({\mathbf{i}^{\prime}}_{d-1},i_d,i^{\prime}_{d+1}))|$.

\item
For some $k$ with $k\neq d$,
$$\begin{cases}
|\overline{M}(S)^b(\mathbf{L}_d(\mathbf{i}_{d+1}))|=|\overline{M}(S)^b(\mathbf{L}_d(\mathbf{i}^{\prime}_{d-1},i_d,i^{\prime}_{d+1}))|\\
(i_{k+1},\ldots,i_{d-1},i_{d+1})=(i^{\prime}_{k+1},\ldots,i^{\prime}_{d-1},i^{\prime}_{d+1}),i_k<i_k^{\prime}
\end{cases}.$$
\end{enumerate}
\end{enumerate}
The definition of $\mathbf{I}$ implies that applying shift operations to the lines $\overline{M}(S)^b(\mathbf{L}_2(\mathbf{i}_{d+1}))$ for $(\mathbf{i}_{d-1},i_{d+1})\in \mathbf{I}$ yields an $(\mathbf{n}_d,q)$-matrix $M^{\prime}$.

To show that the construction of $M^{\prime}$ is valid, we need to verify that the numbers $S(\mathbf{L}_j(\mathbf{i}_{d-1},n_d))$ $(1\leq j\leq d-1$, $\mathbf{i}_{d-1}\in[\mathbf{n}_{d-1}])$ satisfy
$$\overline{S}(\mathbf{L}_j(\mathbf{i}_{d-1},n_d)))\leq S(\mathbf{L}_j(\mathbf{i}_{d-1},n_d))\leq \overline{S}(\mathbf{L}_j(\mathbf{i}_{d-1},1)).$$
If we use the definition of compatibility and the argument used in the proof for the case of $2$-matrices, then we can easily show this.

The line sum array $S(M_1)$ satisfies the assumption of the theorem. Thus induction on the number $d$ guarantees that $M_1$ is constructible.
In addition, the matrix $\overline{M}_{n_d-1}$ is a maximal matrix, hence an $((\mathbf{n}_{d-1},n_d-1),q)$-line sum array $S_{n_d-1}$ defined by
$$\begin{cases}
S_{n_d-1}(\mathbf{L}_j(\mathbf{i}_d))=S(\mathbf{L}_j(\mathbf{i}_d))&(1\leq j\leq d-1,\mathbf{i}_d\in[\mathbf{n}_{d-1}]\times[n_d-1])\\
S_{n_d-1}(\mathbf{L}_d(\mathbf{i}_d))=|\overline{M}_{n_d-1}(\mathbf{L}_d(\mathbf{i}_d))|&(\mathbf{i}_d\in[\mathbf{n}_d])
\end{cases}$$
is compatible.
Therefore, by induction on the number $n_d$, we can construct an $((\mathbf{n}_{d-1},n_d-1),q)$-matrix $M_{n_d-1}$ such that the $(\mathbf{n}_d,q)$-matrix $M=[M_{n_d-1},M_1]$ satisfies
$$S(M)=S.$$
As a result $S$ is valid.
This proves the theorem.
\end{pfqmlmatrix}

Let $\mathfrak{S}_d$ be the permutation set of $[d]$.
For each $\sigma\in \mathfrak{S}_d$, we denote $$\sigma(\mathbf{i}_d)=(i_{\sigma(1)},i_{\sigma(2)},\ldots,i_{\sigma(d)}).$$
An $(\mathbf{n}_d,q)$-matrix $M$ is called \emph{symmetric} if
$$\begin{cases}
n_1=n_2=\cdots=n_d\\
M(\sigma(\mathbf{i}_d))=M(\mathbf{i}_d)\ (\sigma\in\mathfrak{S}_d,\mathbf{i}_d\in[\mathbf{n}_d])
\end{cases}$$
and an $(\mathbf{n}_d,q)$-line sum array $S$ is \emph{symmetric} if
$$\begin{cases}
n_1=n_2=\cdots=n_d\\
S(\mathbf{L}_j(\sigma(\mathbf{i}_d)))=S(\mathbf{L}_j(\mathbf{i}_d))\ (j\in[d],\sigma\in\mathfrak{S}_d,\mathbf{i}_d\in[\mathbf{n}_d])
\end{cases}.$$
We say that a symmetric $(\mathbf{n}_d,q)$-line sum array $S$ is \emph{valid} if there is a symmetric $(\mathbf{n}_d,q)$-matrix $M$ satisfying
$$S(M)=S.$$
Note that a matrix $M$ of size $n_1\times n_2$ is symmetric if
$$n_1=n_2\ \text{and}\ M(i_1,i_2)=M(i_2,i_1).$$
Therefore symmetric multidimensional matrices are generalizations of symmetric $2$-matrices.

Brualdi et al. \cite{brualdi2-1} proved that a criterion for existence of a $q$-ary $2$-matrix with a given line sum array implies a criterion for the case of symmetric $2$-matrices.
Generalizing this result to multidimensional matrices, we  establish a criterion for a $q$-ary symmetric line sum array to be valid.

\begin{theorem}\label{qlsymatrix}
Let $S$ be a symmetric  $q$-ary line sum array.
Then $S$ is valid if and only if $S$ is compatible.
\end{theorem}

\begin{pf}
Let $S$ be a compatible symmetric $(\mathbf{n}_d,q)$-line sum array.
We construct a symmetric $(\mathbf{n}_d,q)$-matrix $M$ with $S(M)=S$ by using double induction on the numbers $d$ and $n_d$.
If $d=2$, then we can easily show that $S$ is valid by applying the proof of Theorem \ref{qlmatrix}.
Thus we assume that $d\ge 3$.
In addition, if $n_d=1$ then $S$ is trivially.
Thus we further assume that $n_d\ge 2$.

We construct an $(\mathbf{n}_d,q)$-matrix $M^{\prime}$ defined as follows:
\begin{enumerate}[1.]
\item
$|M^{\prime}(\mathbf{L}_d(\mathbf{i}_d))|=S(\mathbf{L}_d(\mathbf{i}_d))$  $(\mathbf{i}_d\in[\mathbf{n}_d])$.

\item
$|M^{\prime}(\mathbf{L}_j(\mathbf{i}_{d-1},n_d))|=S(\mathbf{L}_j(\mathbf{i}_{d-1},n_d))$ $(j\in[d-1],\mathbf{i}_{d-1}\in[\mathbf{n}_{d-1}])$.

\item
${M^{\prime}}^b(\sigma(\mathbf{i}_{d-1},n_d),i_{d+1})=1$ $(\sigma\in\mathfrak{S}_d)$ if and only if $(\mathbf{i}_{d-1},i_{d+1})\in \mathbf{I}$.

\item
Let $[(\mathbf{n-1})_d]=\prod_{j=1}^d[n_j-1]$. The $((\mathbf{n-1})_d,q)$-submatrix $\overline{M}_{n_d-1}$ of $M^{\prime}$ defined by
$$\overline{M}_{n_d-1}(\mathbf{i}_d)=M^{\prime}(\mathbf{i}_d)\ (\mathbf{i}_d\in[(\mathbf{n-1})_d])$$ is a maximal matrix.
\end{enumerate}
For this construction,  we define $\mathbf{I}$ to be the set of $(\mathbf{i}_{d-1},i_{d+1})\in[\mathbf{n}_{d-1}]\times[q-1]$ that satisfies the following conditions:
\begin{enumerate}[1.]
\item
For each $1\leq j\leq d-1$ and $\mathbf{i}_{d-1}\in[\mathbf{n}_{d-1}]$, the number of $(\mathbf{i}_{d-1},i_{d+1})$ in $\mathbf{I}$ is $S(\mathbf{L}_j(\mathbf{i}_{d-1},n_d))$.

\item
Any two $(\mathbf{i}_{d-1},i_{d+1})\in\mathbf{I}$ and $(\mathbf{i}^{\prime}_{d-1},i^{\prime}_{d+1})\in[\mathbf{n}_{d-1}]\times [q-1]\setminus \mathbf{I}$ satisfy one of the followings:
\begin{enumerate}[(1)]
\item
$|\overline{M}(S)^b(\mathbf{L}_d(\mathbf{i}_{d+1}))|>|\overline{M}(S)^b(\mathbf{L}_d({\mathbf{i}^{\prime}}_{d-1},i_d,i^{\prime}_{d+1}))|$.

\item
For some $k$ with $k\neq d$,
$$\begin{cases}
|\overline{M}(S)^b(\mathbf{L}_d(\mathbf{i}_{d+1}))|=|\overline{M}(S)^b(\mathbf{L}_d({\mathbf{i}^{\prime}}_{d-1},i_d,i^{\prime}_{d+1}))|\\
(i_{k+1},\ldots,i_{d-1},i_{d+1})=(i^{\prime}_{k+1},\ldots,i^{\prime}_{d-1},i^{\prime}_{d+1}),i_k<i_k^{\prime}
\end{cases}.$$
\end{enumerate}
\end{enumerate}

The $(\mathbf{n}_{d-1},q)$-matrix $M_1$ defined by $$M(\mathbf{i}_{d-1})=|\{i_{d+1}\mid(\mathbf{i}_{d-1},i_{d+1})\in\mathbf{I}\}|$$
for $\mathbf{i}_{d-1}\in[\mathbf{n}_{d-1}]$ is a symmetric matrix whose line sum array satisfies the assumption of the theorem.
Thus induction on the number $d$ ensures that $M_1$ is constructible.
In addition, if we define $S_{n_d-1}$ to be a symmetric $((\mathbf{n-1})_d,q)$-line sum array defined by
$$\begin{cases}
S_{n_d-1}(\mathbf{L}_j(\sigma(\mathbf{i}_d)))=S_{n_d-1}(\mathbf{L}_j(\mathbf{i}_d))&
(j\in[d],\sigma\in \mathfrak{S}_d,\mathbf{i}_d\in[(\mathbf{n-1})_d])\\\
 S_{n_d-1}(\mathbf{L}_d(\mathbf{i}_d)=|\overline{M}_{n_d-1}(\mathbf{L}_d(\mathbf{i}_d))|&( \mathbf{i}_d\in[(\mathbf{n-1})_d])
\end{cases},$$
then, by induction on the number $n_d$, we can construct a symmetric $((\mathbf{n-1})_d,q)$-matrix $M_{n_d-1}$ such that
$$S(M_{n_d-1})=S_{n_d-1}.$$
Therefore the matrix $M$ defined by
$$\begin{cases}
M(\mathbf{i}_d)=M_{n_d-1}(\mathbf{i}_d)&(\mathbf{i}_d\in[(\mathbf{n-1})_d])\\
M(\sigma(\mathbf{i}_{d-1},n_d))=M_1(\mathbf{i}_{d-1})&(\sigma\in\mathfrak{S}_d,\mathbf{i}_{d-1}\in[\mathbf{n}_{d-1}])
\end{cases}$$
is a symmetric $(\mathbf{n}_{d},q)$-matrix such that
$$S(M)=S.$$
This proves the theorem.
\end{pf}

\section*{Classification Codes}
05B20, 15B36

\end{document}